\documentclass[12pt]{amsart}

\usepackage{amscd,amssymb,fullpage}
\usepackage{graphicx}
\usepackage[mathscr]{eucal}

\newcommand{\C}{\mathbb{C}}

\newcommand\of{{\overline{f}}}

\newcommand\oz{{\overline{z}}}
\newcommand\ow{{\overline{w}}}
\newcommand\os{{\overline{s}}}
\newcommand{\Teich}{\textit{Teich}}
\newcommand{\QF}{\textit{QF}}
\newcommand{\M}{\mathcal{M}}
\renewcommand{\S}{\mathcal{S}}
\renewcommand{\H}{\mathbb{H}}
\newcommand{\p}{\partial}
\newcommand{\dbar}{\overline{\partial}}
\DeclareMathOperator{\Spec}{Spec}
\let\Im=\undefined
\DeclareMathOperator{\Im}{Im}

\newcommand{\detp}{\operatorname{det}'}
\DeclareMathOperator{\Diff}{\textit{Diff}}
\DeclareMathOperator{\horizontal}{\mathsf{v}}
\DeclareMathOperator{\vertical}{\mathsf{w}}

\numberwithin{equation}{section}

\theoremstyle{plain}

\newtheorem{theorem}[subsection]{Theorem}
\newtheorem{corollary}[subsection]{Corollary}
\newtheorem{lemma}[subsection]{Lemma}
\newtheorem{proposition}[subsection]{Proposition}

\theoremstyle{remark}

\newtheorem*{example*}{Example}

\newtheorem*{remark*}{Remark}

\begin{document}
\title{Holomorphic extensions of determinants of Laplacians}

\author{Young-Heon Kim}

\address{Department of Mathematics, Northwestern University, Evanston,
  IL 60208, USA}

\email{huns@math.northwestern.edu}

\thanks{The author thanks Ezra Getzler for his advice while this work
  was carried out. He is grateful to Curtis McMullen, Peter Sarnak and
  Jared Wunsch for helpful discussions. This research was partially
  supported by the NSF under grant DMS-9704320}


\begin{abstract}
  The Teichm\"uller space $\Teich(S)$ of a surface $S$ in genus $g>1$
  is a real submanifold of the quasifuchsian space $\QF(S)$. We show
  that the determinant of the Laplacian $\detp(\Delta)$ on $\Teich(S)$
  has a unique holomorphic extension to $\QF(S)$.
\end{abstract}

\maketitle

\section{Introduction} \label{S:setting}

Given a compact Riemannian manifold $M$, the Laplace-Beltrami operator
$\Delta$ on functions on $M$ is an elliptic operator with discrete
spectrum
\begin{equation*}
  \lambda_0 = 0 < \lambda_1 \leq \lambda_2 \leq \cdots \leq
  \lambda_k \leq \cdots \rightarrow \infty .
\end{equation*}
The determinant of the operator $\Delta$ may be defined formally as
the product of the nonzero eigenvalues of $\Delta$. A regularization
$\detp(\Delta)$ of this product was defined by Ray and Singer
\cite{RS1}, \cite{RS2}, using the zeta function of $\Delta$.

The determinant $\detp(\Delta)$ has played an important role in such
areas of mathematics as mathematical physics, differential geometry,
algebraic geometry and number theory. In particular, it plays a
central role in the theory of determinant line bundles, initiated by
Quillen (\cite{Q}) and further developed by Bismut and Freed
\cite{BF1}, \cite{BF2}, and by Bismut, Gillet and Soul\'e \cite{BGS1},
\cite{BGS2}, \cite{BGS3}.

In a series of papers \cite{OPS1}, \cite{OPS2} (see also \cite{S2}),
Osgood, Phillips and Sarnak studied $-\log\detp(\Delta)$ as a
``height'' function on the space of metrics on a compact orientable
smooth surface $S$ of genus $g$. For $g>1$, they showed that when
restricted to a given conformal class of metrics on $S$, it attains
its minimum at the unique hyperbolic metric in this conformal class,
and has no other critical points. Thus, to find Riemannian metrics on
$S$ which are extremal, in the sense that they minimize
$-\log\detp(\Delta)$, it suffices to consider its restriction to the
moduli space $\M_g$ of hyperbolic metrics on a Riemann surface $S$ of
genus $g$. Osgood, Phillips and Sarnak showed that this restriction is
a proper function.

The universal cover of the orbifold $\M_g$, with covering group the
mapping class group $\Gamma_g$, is the Teichm\"uller space
$\Teich(S)$. The function $-\log\detp(\Delta)$ lifts to a function on
the Teichm\"uller space $\Teich(S)$ invariant under $\Gamma_g$.  In
this paper, we are interested the function theoretic properties of
$\log\detp(\Delta)$ on $\Teich(S)$. Before stating the main theorem,
consider the special case of genus $1$.

\begin{example*}[\cite{RS2} or \cite{S1}, p.\ 33, (A.1.7)]
  For $z \in \H$, let $T_z$ be the flat torus obtained by the
  lattice of $\C$ generated by $1$ and $z$. Then the
  determinant of Laplacian of this flat torus is
  \begin{equation*}
  \log \detp(\Delta)(z) = \log (2\pi y^{1/2} |\eta (z) |^2 )
  \end{equation*}
  where $\eta(z)= q^{1/24} \prod_{n=1}^{\infty} (1- q^{n})$ for
  $q=e^{2\pi i z}$ is the eta function; this is a modular form of
  weight $1/2$.

  The manifold $\H$ has a complexification $\H \times
  \overline{\H}$, and the function $\log\detp(\Delta)(z)$ on
  the diagonal $\{w=\overline{z}\}$ has a holomorphic extension to
  $\H \times \overline{\H}$, namely,
  \begin{equation*}
  \log ( - \pi i(z-w)^{1/2} \, \eta(z) \, \overline{\eta(\overline{w})} ) .
  \end{equation*}
\end{example*}

In this paper, we show that even in higher genus $g>1$, the function
$\log \detp(\Delta)$ has a unique holomorphic extension. In higher
genus, the objects corresponding to $\H$ and
$\H\times\overline{\H}$ are the Teichm\"uller space
$\Teich(S)$ and the quasifuchsian space
$\QF(S)=\Teich(S)\times\Teich(\overline{S})$ respectively. The complex
structure on quasifuchsian space $\QF(S)$ is induced by the complex
structure of $\Teich(S)$, and the real analytic manifold $\Teich(S)$
imbeds as the diagonal in $\QF(S)$.

The following theorem is the main result of this paper. (For the
precise statement, see Theorem~\ref{holo det in quasifuchsian}.)
\begin{theorem} \label{main}
  The function $\log\detp(\Delta)$ on $\Teich(S)$ has a unique
  holomorphic extension to the quasifuchsian space $\QF(S)$.
\end{theorem}

In the proof of Theorem~\ref{main}, we use the Belavin-Knizhnik
formula (see Theorem~\ref{curvature}), proved by Wolpert \cite{W3} and
by Zograf and Takhtajan \cite{ZT} and the holomorphic extension of the
Weil-Petersson form constructed by I.~Platis~\cite{P} (see
Theorem~\ref{Platis}).

The following is a key step in the proof of Theorem~\ref{main} and  may be
of independent interest.
\begin{theorem}
  Let $V$ and $W$ be domains in the complex space $\C^n$ diffeomorphic
  to the open unit ball. Consider $V \times W\subset\C^n\times\C^n$,
  with holomorphic coordinates $(z, w)$, and let $\p_z = dz^i \,
  \p_{z_i}$ and $\p_w = dw^j\, \p_{w_j}$. Suppose
  $\Omega$ is a holomorphic closed 2-form on $V\times W$ which is
  locally written as
  \begin{equation*}
  \Omega = \sum_{i,j} \Omega_{ij} dz^i \wedge dw^j .
  \end{equation*}
  Then there is a holomorphic function $q$ on $V \times W$ such
  that $\p_z \p_w q = \Omega$.
\end{theorem}

This theorem implies that there is a holomorphic function on $\QF(S)$
whose restriction on $\Teich(S)$ (the diagonal in $\Teich(S) \times
\Teich(\overline{S})$) is a K\"ahler potential for the Weil-Petersson
form $\omega_{WP}$. This suggests that quasifuchsian space is a useful
tool in gaining a better understanding of function theory of the
Teichm\"uller space. This idea is due to McMullen
\cite{M}, who used quasifuchsian space to study the geometry of
the Weil-Petersson metric on the Teichm\"uller space.

In a sequel to this paper, we will construct a holomorphic family of
differential operators in a neighborhood of the diagonal in
quasifuchsian space, which equal the Laplacian along the diagonal, and
such that the determinant of this family equals the holomorphic
extension of $\log\detp(\Delta)$ constructed in this paper.

The asymptotic behavior of $\log\detp(\Delta)$ near the boundary of
Teichm\"uller space is important in both geometry and physics and was
studied in \cite{W4} and \cite{BB}. It would be interesting to
understand the asymptotic behavior of the holomorphic extension of
$\log\detp(\Delta)$ near the boundary of the quasifuchsian space. We
hope to address this in the future.

\subsection*{Plan of the paper}
In Section 2, we review the facts that we need on Teichm\"uller spaces
and quasifuchsian spaces, including the Belavin-Knizhnik formula and
Platis's theorem. We prove Theorem~\ref{potential} in Section 3.  In
Section 4, we complete the proof of Theorem~\ref{main}.

\section{Preliminaries}

\subsection*{Determinants of Laplacians}

Let $\Delta $ be the Laplace-Beltrami operator on functions on a
compact Riemannian manifold $M$. Let
\begin{equation} \label{zeta}
  \zeta_\Delta(s)=  \sum_{\lambda \in \Spec(\Delta)\setminus\{0\}}
  \lambda^{-s}
\end{equation}
be the zeta-function of $\Delta$. The determinant $\detp(\Delta)$ is
defined (see \cite{RS1}) as
\begin{equation} \label{det'}
  - \log \detp(\Delta) = \frac{d\zeta_\Delta(0)}{ds} .
\end{equation}
Since $M$ is compact and $\Delta$ is elliptic and self-adjoint, the
nonzero spectrum of $\Delta$ is positive and discrete. Moreover, the
sum in Example~\ref{zeta} is absolutely convergent for $\Re s$
sufficiently large, and has a meromorphic extension to the whole
complex plane, with possible poles only at $\{ -1, -2, -3, \dots \}$.
Thus, there is no difficulty in taking the derivative at $s=0$ in
\eqref{det'}.

\subsection*{Teichm\"uller spaces}

A general reference for this section is \cite{IT}.

Let $S$ be an oriented closed surface with genus $g > 1$. The
Teichm\"uller space $\Teich(S)$ of $S$ is the space of isotopy classes
of hyperbolic Riemannian metrics on $S$, that is, metrics with
Gaussian curvature $-1$. Two Riemannian metrics $m_1$ and $m_2$ on $S$
are said to be in the same isotopy class if there is an isotopy $\phi$
of $S$, i.e. a diffeomorphism obtained by a flow of a vector field on
$S$, such that $\phi^* m_1 = m_2$. On a surface, there is one-to-one
correspondence between complex structures and hyperbolic Riemannian
metrics, i.e.\ each complex structure on $S$ has unique hyperbolic
metric which is a K\"ahler metric to the complex structure and each
hyperbolic Riemannian metric on $S$ has canonical complex structure
such that the metric becomes K\"ahler. From this correspondence we see
that $\Teich(S)$ is also the space of isotopy classes of complex
structures on $S$.

The set of equivalence classes of hyperbolic metrics (or equivalently
complex structures) under orientation preserving diffeomorphisms on
$S$ forms the moduli space $\M_g$ of compact Riemann surfaces of genus
$g$.

Denote the group of orientation preserving diffeomorphisms on $S$ by
$\Diff^+(S)$, and the group of isotopies by $\Diff_0(S)$. The mapping
class group
\begin{equation*}
\Gamma_g = \Diff^+(S) / \Diff_0 (S)
\end{equation*}
is a discrete group which acts properly discontinuously on
$\Teich(S)$. Thus $\Teich(S)$ is almost a covering space of $\M_g$, with
covering transformation group $\Gamma_g$:
\begin{equation*}
\begin{CD}
  \Gamma_g @>>> \Teich(S)  \\
  @. @VVV \\
  @. \M_g = @. \Gamma_g \backslash \Teich(S)
\end{CD}
\end{equation*}
The only caveat is that the action of $\Gamma_g$ is not free, i.e.\
there are points in $\Teich(S)$ which are fixed under some finite
subgroups of $\Gamma_g$. These points descend to $\M_g$ as orbifold
singularities.

Fixing a hyperbolic metric on $S$, we may decompose $S$ into $2g-2$
pairs of pants, separated by closed geodesics
$\gamma_1,\dots,\gamma_{3g-3}$.

A hyperbolic pair of pants is determined up to isometry by the lengths
of its boundary geodesics. Given the combinatorial pants decomposition
of $S$, we get a hyperbolic metric by specifying the lengths $l_i$
($l_i >0 $) of the geodesics $\gamma_i$ and the angle $\theta_i$ by
which they are twisted along $\gamma_i$ before gluing. Let
$\tau_i=l_i\theta_i/2\pi$, $i=1,\dots,3g-3$. Then the system of
variables
\begin{equation*}
(l_1 , \dots , l_{3g-3}, \tau_1 , \dots , \tau_{3g-3} )
\end{equation*}
is a real analytic coordinate system on $\Teich(S)$, called the
\emph{Fenchel-Nielsen coordinates} of $\Teich(S)$. This coordinate
system gives a diffeomorphism
\begin{equation*}
\Teich(S) \approx \mathbb{R}_+ ^{3g-3} \times \mathbb{R}^{3g-3} \, .
\end{equation*}

There is a a natural symplectic form $\omega_{WP}$ on $\Teich(S)$,
called the Weil-Petersson form. By a theorem of Wolpert (\cite{W1},
\cite{W2}; see also \cite{IT}), this form is given in Fenchel-Nielsen
coordinates by the formula
\begin{equation} \label{WP}
  \omega_{WP} = \sum_{i=1}^{3g-3} dl_i \wedge \, d\tau_i .
\end{equation}

The Teichm\"uller space $\Teich(S)$ has a natural complex structure,
for which $\omega_{WP}$ is a K\"ahler form. The following theorem is
well known. (See, for example, \cite{A}.)
\begin{theorem} \label{Teich as domain}
  For a closed surface $S$ with genus $g>1$, $\Teich(S)$ is
  biholomorphic to a bounded open contractible domain in $\C^{3g -3}$.
\end{theorem}

\begin{corollary} \label{global coordinate}
  There are global holomorphic coordinates $z=(z^1,\dots,z^{3g-3})$ on
  $\Teich(S)$.
\end{corollary}


\subsection*{Quasifuchsian spaces}
While Teichm\"uller space is a space of Riemann surfaces, the
quasifuchsian space defined by Lipman Bers (See \cite{B}) is a space
of pairs of Riemann surfaces. The quasifuchsian space $\QF(S)$ of the
surface $S$ may simply be defined as
\begin{equation*}
\QF(S) = \Teich(S) \times \Teich(\overline{S}) .
\end{equation*}
Here, $\overline{S}$ denotes the surface $S$ with the opposite
orientation.

The complex conjugate $\overline{X}$ of a Riemann surface $X$ is
defined by the following diagram:
\begin{equation} \label{conjugation}
  \begin{array}{clcl}
    \H & \longrightarrow & \overline{\H} & \\
    \downarrow & & \downarrow & \\
    X & \dashrightarrow & \overline{X} & \\
  \end{array}
\end{equation}
The upper arrow is complex conjugation, and the vertical arrows are
the universal coverings given by the uniformization theorem for
Riemann surfaces. There is a canonical map from $\Teich(S)$ to
$\Teich(\overline{S})$ defined by sending a Riemann surface
$X\in\Teich(S)$ to its complex conjugate
$\overline{X}\in\Teich(\overline{S})$.

Let $\overline{\Teich(S)}$ be the complex conjugate of $\Teich(S)$.
\begin{proposition} \label{Teich(bar-S)}
  As complex manifolds, $\Teich(\overline{S}) \cong
  \overline{\Teich(S)}$.
\end{proposition}
\begin{proof}
  The complex structure of $\Teich(S)$ is induced by the complex
  structure of the space of Beltrami differentials on $S$. Fix a
  complex structure on $S$. Beltrami differentials are $(-1,1)$-forms
  on $S$ with $L_\infty$ norm less than $1$; locally, $\mu(z)\,d\bar
  z/dz $ with $|\mu|<1$. (See \cite{A} or \cite{AB} for details.)

  If $z$ is a local coordinate on $S$, then $w=\overline{z}$ is a
  local coordinate on $\overline{S}$. Now Beltrami differentials on
  $\overline{S}$ are locally of the form $\overline{\mu(\overline{w})}
  d\overline{w}/dw$. From this local expression it is clear that the
  complex structure on $\Teich(\overline{S})$ is the same as the one
  on $\overline{\Teich(S)}$.
\end{proof}

The diagonal map
\begin{equation*}
  \Teich(S) \hookrightarrow \Teich(S) \times \overline{\Teich(S)}
\end{equation*}
sending $X\in\Teich(S)$ to $(X,\overline{X})$ embeds $\Teich(S)$ as
a totally real submanifold into $\QF(S)$. The action of $\Gamma_g$ on
$\Teich(S)$ extends to $\QF(S)$ by the diagonal action: for $\rho \in
\Gamma_g$ and $(X,Y) \in \QF(S)= \Teich(S)\times\overline{\Teich(S)}$,
\begin{equation*}
  \rho\cdot(X, Y) = (\rho\cdot X , \rho \cdot Y) .
\end{equation*}

By Corollary~\ref{global coordinate}, $\QF(S)=\Teich(S) \times
\overline{\Teich(S)}$ has global holomorphic coordinates
\begin{equation*}
  (z^1 , \dots , z^{3g-3} , w^1 , \dots , w^{3g-3}) .
\end{equation*}
We abbreviate this coordinate system to $(z,w)$.

\subsection*{Holomorphic extension of Weil-Petersson form}
The following result is due to Platis (\cite{P}, Theorems 6 and 8).
\begin{theorem} \label{Platis}
  The differential form $i\omega_{WP}$ on the Teichm\"uller space
  $\Teich(S)$ has an extension $\Omega$ to the quasifuchsian space
  $\QF(S)$ which is a holomorphic non-degenerate closed $(2,0)$-form
  whose restriction to the diagonal
  $\Teich(S)\subset\QF(S)\cong\Teich(S)\times\overline{\Teich(S)}$ is
  $i\omega_{WP}$.
\end{theorem}

\begin{lemma} \label{identity theorem}
  Let $U\subset\C^n$ be a connected complex domain, and let $\phi$ be
  a holomorphic function on $U\times\overline{U}$ whose restriction to
  the diagonal $U\subset U\times\overline{U}$ vanishes. Then $\phi$
  vanishes on all of $U\times\overline{U}$.
\end{lemma}
\begin{proof}
  It suffices to show that $\phi$ vanishes near a point $(z,\oz)$ on
  the diagonal. Choose holomorphic coordinates $s$ on $U$ which vanish
  at $z$, and let $(s,t)$ be the corresponding holomorphic coordinates
  on $U\times\overline{U}$. Consider the Taylor series expansion
  \begin{equation*}
    \phi(s,t) = \sum_{\alpha,\beta} a_{\alpha,\beta} \, s^\alpha t^\beta
    .
  \end{equation*}
  Along the diagonal, where $\os=t$, we have
  \begin{equation*}
    0 = \phi(s,\os) = \sum_{\alpha,\beta} a_{\alpha,\beta} \, s^\alpha
    \os^\beta .
  \end{equation*}
  It follows that $a_{\alpha,\beta}=0$ for all $\alpha$ and
  $\beta$, hence $\phi(s,t)$ vanishes in a neighborhood of $(z,\oz)$.
\end{proof}

\begin{proposition} \label{form of Omega}
  In terms of the holomorphic coordinate system
  \begin{equation*}
    (z,w)=(z^1, \dots, z^{3g-3}, w^1, \dots , w^{3g-3} )
  \end{equation*}
  on $\Teich(S) \times \overline{\Teich(S)}$, the $2$-form $\Omega$ of
  Theorem~\ref{Platis} may be written locally as
  \begin{equation*}
  \Omega = \sum_{i,j} \Omega_{ij} \, dz^i \wedge dw^j .
  \end{equation*}
\end{proposition}
\begin{proof}
  Since $\Omega$ is $(2,0)$ form, we may write
  \begin{equation*}
    \Omega = \sum_{i,j} \bigl( A_{ij} \, dz^i \wedge dz^j + B_{ij} \,
    dz^i \wedge dw^j + C_{ij} \, dw^i \wedge dw^j \bigr) .
  \end{equation*}
  Restricted to the diagonal $\{w=\oz\}$,
  \begin{equation*}
    i\omega_{WP} = \Omega|_{w=\oz} = \sum_{i,j}
    \bigl( A_{ij}|_{w=\oz} \, dz^i \wedge dz^j + B_{ij}|_{w=\oz} \, dz^i
    \wedge d\oz^j + C_{ij}|_{w=\oz} \, d\oz^i \wedge d\oz^j \bigr) .
  \end{equation*}
  Since $\omega_{WP}$ is $(1,1)$-form on $\Teich(S)$, we see that
  $A_{ij}$ and $C_{ij}$ vanish on the diagonal $\{w=\overline{z}\}$.
  Since $\Omega$ is holomorphic, so are the functions $A_{ij}$,
  $B_{ij}$, and $C_{ij}$. Applying Lemma~\ref{identity theorem}, we
  see that $A_{ij}$ and $C_{ij}$ vanish.
\end{proof}

\subsection*{The Laplacian on hyperbolic surfaces and the
  Belavin-Knizhnik formula}

Let $X$ be a compact hyperbolic surface of genus $g>1$, and let
$\Delta$ be the Laplacian on scalar functions on $X$. On the universal
cover $\H$ of $X$, the pull-back of $\Delta$ by the covering
map may be written as
\begin{equation*}
\Delta = (z-\overline{z})^2 \frac{\p^2}{\p z \p{\overline z}} .
\end{equation*}

The \emph{Siegel upper half space} $\S_g$ is the space of complex
symmetric matrices in $\C^{g\times g}$ with positive definite
imaginary part. The period matrix $\tau$ is a holomorphic map from
$\Teich(S)$ to $\S_g$.

We will use the Belavin-Knizhnik formula, proved by Wolpert and by
Zograf and Takhtajan. (See \cite{W3} and \cite{ZT}.)  We only need the
following special case of this theorem (\cite{ZT}, Theorem 2).
\begin{theorem} \label{curvature}
  In $\Teich(S)$,
  \begin{equation*}
    \p \dbar \, \log \left( \frac{\detp(\Delta)}{\det(\Im\tau)}
    \right) = iC_g \, \omega_{WP} ,
  \end{equation*}
  where $\Im\tau$ is the imaginary part of the period matrix $\tau$
  and $C_g$ is a constant depending only on the genus $g$. The
  differential operator $\p\dbar$ comes from the complex structure on
  $\Teich(S)$.
\end{theorem}

This formula and the result of the next section together with the
theorem of Platis are the key ingredients in the construction of the
holomorphic extension of $\log\detp(\Delta)$.

\section{An extended K\"ahler potential} \label{S:almas}

The proof of the following elementary theorem occupies the whole of
this section.
\begin{theorem} \label{potential}
  Let $V$ and $W$ be domains in the complex space $\C^n$ which are
  diffeomorphic to the open unit ball, and let
  $(z,w)=(z^1,\dots,z^n,w^1,\dots,w^n)$ be holomorphic coordinates on
  the domain $V\times W\subset\C^n\times\C^n$. Suppose $\Omega$ is a
  holomorphic closed 2-form on $V\times W$ of the form
  \begin{equation*}
  \Omega = \sum_{i,j} \Omega_{ij} \, dz^i \wedge dw^j .
  \end{equation*}
  Then there exists a holomorphic function $q$ on $V\times W$ such
  that
  \begin{equation*}
  \p_z \p_w q = \Omega ,
  \end{equation*}
  where $\p_z=\sum_idz^i\,p_{z^i}$ and $\p_w=\sum_jdw^j\,\p_{w^j}$.
\end{theorem}
\begin{proof}
  Choose smooth polar coordinates on $V$ and $W$, and denote the
  centers of these coordinate systems by $z_0$ and $w_0$ respectively.
  Denote the radial line in polar coordinates from $z_0$ to the point
  $z\in V$ by $\horizontal(z)$; similarly, denote the radial line in
  polar coordinates from $w_0$ to the point $w\in W$ by
  $\vertical(w)$. More generally, if $c$ is a smooth chain in $V$, let
  $\horizontal(c)$ denote the cone on $c$ with vertex $z_0$, and
  similarly if $c$ is a smooth chain in $W$, let $\vertical(c)$ denote
  the cone on $c$ with vertex $w_0$.

  Define $q(z,w)$ by the formula
  \begin{equation*}
    q(z,w) = \int_{\horizontal(z) \times \vertical(w)} \Omega .
  \end{equation*}
  Since the cycle $\horizontal(z) \times \vertical(w)$ varies smoothly
  as $(z,w)$ varies, the function $q(z,w)$ is smooth. Observe that $q$
  is unchanged by isotopies of the coordinate systems on $V$ and $W$
  which fix the centers $z_0$ and $w_0$, and that $q$ vanishes on
  $V\times\{w_0\}$ and on $\{z_0\}\times W$.

  If $c$ is a differentiable curve in $W$ parametrized by the interval
  $[0,t]$, we have by Stokes's theorem
  \begin{equation*}
    q(z,c(t)) - q(z,c(0)) = \int_{\horizontal(z)\times c} \Omega +
    \int_{\{z\}\times \vertical(c)} \Omega
    - \int_{\{z_0\} \times\vertical(c)} \Omega -
    \int_{\horizontal(z) \times \vertical(c)} d\Omega .
  \end{equation*}
  The second and third terms on the right-hand side vanish, since
  $\Omega$ vanishes when restricted to the $2$-simplex
  $\{z\}\times\vertical(c)$, and the last term vanishes since
  $d\Omega=0$. Taking the limit $t\to0$, we see that
  \begin{equation} \label{d_w} \iota(0,c'(0)) dq(z,c(0)) = -
    \int_{\horizontal(z)\times c(0)} \iota(0,c'(0))\Omega .
  \end{equation}
  Since $\Omega$ is holomorphic along $\{z\} \times W$, it follows
  that $q$ is holomorphic along $\{z\}\times W$ as well. A similar
  argument shows that $q$ is holomorphic along $V\times\{w\}$;
  combining these two calculations, we see that $q$ is holomorphic on
  $V\times W$.

  We now calculate $\p_w\p_zq$. By \eqref{d_w},
  \begin{equation*}
    \p_w q(z,w)  = -\sum_{i=1}^n dw^i  \int_{\horizontal(z)\times\{w\}}
    \iota(\p_{w^i}) \Omega .
  \end{equation*}
  If $c$ is a differentiable curve in $V$, parametrized by the
  interval $[0,t]$, we have by Stokes's theorem
  \begin{equation} \label{d_z}
    (\p_wq)(c(t),w) - (\p_wq)(c(0),w) = \sum_{i=1}^n dw^i \left( -
    \int_{c\times\{w\}} \iota(\p_{w^i}) \Omega
    + \int_{\horizontal(c)\times\{w\}} d\iota(\p_{w^i})\Omega \right) .
  \end{equation}
  The second term on the right-hand side vanishes. Indeed,
  \begin{align*}
    d \iota(\p_{w^i}) \Omega &= - \p_{z^k} \Omega_{ji} dz^k \wedge
    dz^j - \p_{w^k} \Omega_{ji} dw^k \wedge dz^j \\
    &= - \sum_{j<k} \bigl( \p_{z^k} \Omega_{ji} - \p_{z^j} \Omega_{ki}
    \bigr) \, dz^k \wedge dz^j - \p_{w^k} \Omega_{ji} dw^k \wedge dz^j \\
    &= - \p_{w^k} \Omega_{ji} dw^k \wedge dz^j .
  \end{align*}
  Restricting to $\horizontal(c)\times\{w\}$, this differential form
  vanishes.

  Taking $t\to0$ in \eqref{d_z}, we see that
  \begin{equation*}
    \iota(c'(0),0)(\p_wq)(c(0),w) = - \sum_{i=1}^n dw^i \,
    \iota(c'(0),0) \iota(\p_{w^i}) \Omega(c(0),w) ,
  \end{equation*}
  or in other words, $\p_z\p_wq=\Omega$.
\end{proof}

\section{Holomorphic extension of $\log\detp\Delta$} \label{SS:2mflds}

From Proposition~\ref{form of Omega}, we know that the holomorphic
2-form $\Omega$ of Theorem~\ref{Platis} satisfies the hypotheses of
Theorem~\ref{potential}. Restricted to the diagonal
$\Teich(S)=\{w=\overline{z}\}\subset\QF(S)$, the differential equation
in Theorem~\ref{potential} for the holomorphic function $q$ on
$\QF(S)$ becomes
\begin{equation*}
  \p \dbar q = i \omega_{WP} ,
\end{equation*}
where $i\omega_{WP}$ is the restriction of $\Omega$ to the diagonal.
Thus, Theorem~\ref{potential} gives a method of constructing a
K\"ahler potential for the K\"ahler form $i\omega_{WP}$ on the
Teichm\"uller space, using the fact that it has a holomorphic
extension to quasifuchsian space.

\begin{example*}  (See p214 in \cite{IT})
  When $S$ has genus $1$, the Teichm\"uller space $\Teich(S)$ may be
  identified with the upper half plane $\H$, and
  \begin{equation*}
  \omega_{WP} = -i (z-\overline{z})^{-2} \, dz \wedge d\overline{z} \  .
  \end{equation*}
One easily finds the K\"ahler potential $q(z)=\log(z-\oz)$. The method
  used in the proof of Theorem~\ref{potential}, applied to the 2-form
  $\Omega = (z-w)^{-2} \, dz \wedge dw$, yields the holomorphic
  function
  \begin{equation*}
    q(z,w) = \log(z-w) - \log(z_0-w) - \log(z-\ow_0) + \log(z_0-\ow_0))
  \end{equation*}
  on the quasifuchsian space $\H\times\overline{\H}$.
\end{example*}

\begin{lemma} \label{real potential}
  There is a holomorphic function $\tilde{q}(z,w)$ on the
  quasifuchsian space
  $\QF(S)\cong\Teich(S)\times\overline{\Teich(S)}$, whose restriction
  to the diagonal $\Teich(S)=\{w=\overline{z} \}$, $\tilde{q}$ is
  real, and such that
  \begin{equation*}
    \p \dbar \tilde{q} = i \omega_{WP} .
  \end{equation*}
\end{lemma}
\begin{proof}
  The function $\overline{q(\ow,\oz)}$ is holomorphic, and
  \begin{equation*}
    \p_z \p_w \bigl( q(z,w) + \overline{q(\ow,\oz)} \bigr) =
    \Omega(z,w) - \overline{\Omega(\ow,\oz)} .
  \end{equation*}
  Restricted on the diagonal $\{w=\oz\}$, we have
  \begin{equation*}
    \p \dbar \bigl( q(z,\oz) + \overline{q(z,\oz)} \bigr) = 2i\omega_{WP} .
  \end{equation*}
  Thus, it suffices to take $\tilde{q}(z,w) = \frac{1}{2} \bigl(
  q(z,w) + \overline{q(\ow,\oz)} \bigr)$.
\end{proof}

\begin{proposition} \label{holo det in Teich} If genus $g>1$, there is
  a holomorphic function $f$ on $\Teich(S)$ such that, in the notation
  of Theorem~\ref{curvature},
  \begin{equation} \label{real log det}
    \log \detp(\Delta) =  C_g \tilde{q} + \log\det(\Im\tau) + f +
    \overline{f} .
  \end{equation}
\end{proposition}
\begin{proof}
  Let $v$ be the real function
  \begin{equation*}
    v = - C_g \tilde{q} + \log \left( \frac{\detp(\Delta)}
      {\det(\Im\tau)} \right)
  \end{equation*}
  on $\Teich(S)$. By Theorem~\ref{curvature} and Lemma~\ref{real
    potential},
  \begin{equation*}
    d (\p v) = 0 .
  \end{equation*}
  Since $\Teich(S)$ is diffeomorphic to an open ball, there is a
  smooth function $h$ on $\Teich(S)$ such that
  \begin{equation*}
    df = \p v .
  \end{equation*}
  It follows that $f$ is holomorphic on $\Teich(S)$, and $v=f+\of$, up
  to a constant which may be absorbed into the definition of $f$.
\end{proof}

This proposition gives rise to the holomorphic extension of
$\log\detp(\Delta)$, since each of the terms on the right-hand side of
\eqref{real log det} has a natural holomorphic extension to $\QF(S)$.
\begin{theorem} \label{holo det in quasifuchsian} There exists a
  unique holomorphic extension of $\log\detp(\Delta)$ to the
  quasifuchsian space $\QF(S)$. In coordinates $(z,w)$ on
  $\QF(S)\cong\Teich(S)\times\overline{\Teich(S)}$, this extension has
  the form
  \begin{equation*}
    \log\detp(\Delta)(z,w) =  C_g \tilde{q}(z,w) +
    \log\det\bigl( (\tau(z)-\overline{\tau(\ow)})/2i \bigr) + f(z) +
    \overline{f(\ow)} .
  \end{equation*}
\end{theorem}
\begin{proof}
  The only term whose extension is not obvious is
  \begin{equation*}
    \log\det(\Im\tau) = \log \det\bigl( ( \tau - \overline{\tau} )/2i
    \bigr) .
  \end{equation*}
  This has the holomorphic extension
  \begin{equation*}
    \log \det\bigl( ( \tau(z) - \overline{\tau(\ow)} )/2i  \bigr) ;
  \end{equation*}
  we need only to observe that the matrix $\tau(z)-\overline{\tau(\ow)}$
  is everywhere invertible on $\QF(S)$.

  The uniqueness of the holomorphic extension of $\log\detp(\Delta)$
  follows from Lemma \ref{identity theorem}.
\end{proof}

\begin{remark*}
  We know that $\log\detp(\Delta)$ on $\Teich(S)$ is invariant under
  the action of the mapping class group. It follows that its
  holomorphic extension to
  $\QF(S)\cong\Teich(S)\times\overline{\Teich(S)}$ is invariant under
  the diagonal action of the mapping class group, because this action
  is holomorphic.
\end{remark*}

\begin{remark*}
  Theorem \ref{holo det in quasifuchsian} implies that
  $\log\detp(\Delta)$ is a real analytic function on $\Teich(S)$, and
  in fact, may be used to give a lower bound for its radius of
  convergence.
\end{remark*}

\bibliographystyle{plain}

\end{document}